\documentclass{amsart}

\usepackage{amsmath,amscd,amssymb,amsthm}
\usepackage{enumerate,mathrsfs}
\usepackage{geometry}
\usepackage{hyperref}
\usepackage{url}
\usepackage[all]{xy}

\newcommand{\bibfilename}{/home/micromath/results/b20140922unibib/unibib}

\numberwithin{equation}{section}

\theoremstyle{plain}
\newtheorem{thm_}[equation]{Theorem}
\newtheorem{lemma_}[equation]{Lemma}
\newtheorem{prop_}[equation]{Proposition}
\newtheorem{cor_}[equation]{Corollary}
\newtheorem{eg_}[equation]{Example}
\newtheorem{con_}[equation]{Conjecture}
\newtheorem*{cons_}{Conjecture}

\theoremstyle{definition}
\newtheorem{thmu_}[equation]{Theorem}
\newtheorem*{thmus_}{Theorem}
\newtheorem{propu_}[equation]{Proposition}
\newtheorem*{propus_}{Proposition}
\newtheorem{coru_}[equation]{Corollary}
\newtheorem{lemu_}[equation]{Lemma}
\newtheorem*{lemus_}{Lemma}
\newtheorem{egu_}[equation]{Example}
\newtheorem*{egus_}{Example}
\newtheorem{def_}[equation]{Definition}
\newtheorem*{defs_}{Definition}
\newtheorem{rk_}[equation]{Remark}

\newcommand{\thm}[1]{\begin{thm_}#1\end{thm_}}

\newcommand{\thmus}[1]{\begin{thmus_}#1\end{thmus_}}
\newcommand{\lemm}[1]{\begin{lemma_}#1\end{lemma_}}

\newcommand{\eg}[1]{\begin{eg_}#1\end{eg_}}

\newcommand{\prop}[1]{\begin{prop_}#1\end{prop_}}

\newcommand{\rk}[1]{\begin{rk_}#1\end{rk_}}

\newcommand{\coru}[1]{\begin{coru_}#1\end{coru_}}

\newcommand{\pf}[1]{\begin{proof}#1\end{proof}}

\DeclareMathOperator{\sgn}{sgn}

\DeclareMathOperator{\Gal}{Gal}

\DeclareMathOperator{\Spec}{Spec}

\newcommand{\ZZ}{\mathbb Z}
\newcommand{\QQ}{\mathbb Q}

\newcommand{\FF}{\mathbb F}
\newcommand{\bA}{\mathbb A}%
\newcommand{\II}{\mathbb I}
\newcommand{\GG}{\mathbf G}

\newcommand{\TT}{\mathbf T}
\newcommand{\XX}{\mathbf X}

\newcommand{\fl}{\mathfrak l}%

\newcommand{\fo}{\mathfrak o}%
\newcommand{\p}{\mathfrak p}

\newcommand{\fr}{\mathfrak r}%

\newcommand{\fL}{\mathfrak L}%
\newcommand{\fP}{\mathfrak P}%

\newcommand{\s}{\sigma}

\newcommand{\tm}{\times}%
\newcommand{\otm}{\otimes}

\newcommand{\lra}{\longrightarrow}

\newcommand{\lmpt}{\longmapsto}

\newcommand{\fracn}[2]{\genfrac{(}{)}{}{}{#1}{#2}}

\newcommand{\eq}[1]{\begin{equation}#1\end{equation}}
\newcommand{\eqn}[1]{\begin{equation*}#1\end{equation*}}

\newcommand{\gan}[1]{\begin{gather*}#1\end{gather*}}
\newcommand{\al}[1]{\begin{align}#1\end{align}}
\newcommand{\aln}[1]{\begin{align*}#1\end{align*}}

\newcommand{\cs}[1]{\begin{cases}#1\end{cases}}

\newcommand{\enmt}[1]{\begin{enumerate}#1\end{enumerate}}

\newcommand{\pin}{{\p_\infty}}
\newcommand{\Pin}{{\fP_\infty}}
\newcommand{\Ep}{{E_\Pin^+}}
\newcommand{\Xp}{{\Xi_\Pin^+}}
\newcommand{\tXp}{{\tilde\Xi_\Pin^+}}
\newcommand{\Kp}{{K_\Pin^+}}
\newcommand{\EP}{{E_\Pin^\tm}}
\newcommand{\XP}{{\Xi_\Pin}}
\newcommand{\KP}{{K_\Pin}}
\DeclareMathOperator{\lc}{lc}

\begin{document}
\title[Diophantine equations defined by binary quadratic forms]{Diophantine
equations  defined by  binary quadratic forms over  rational function fields}
\author[C. Lv]{Chang Lv}
\address{State Key Laboratory of Information Security\\
Institute of Information Engineering\\
Chinese Academy of Sciences\\
Beijing 100093, P.R. China}
\email{lvchang@amss.ac.cn}
\date{\today}
\thanks{This work was supported by
 National Natural Science Foundation of China (Grant No. 11701552).
}
\begin{abstract}
We study the ``imaginary" binary quadratic
 form equations $ax^2+bxy+cy^2+g=0$ over  $k[t]$ in rational function fields,
 showing that a condition with respect to the Artin reciprocity map, is the only
 obstruction to the local-global principle for integral solutions of the equation.
\end{abstract}
\maketitle

\section{Introduction}\label{sec_intro}
Consider the integral solvability of the generalized equation
\eq{\label{eq_qr}
ax^2+bxy+cy^2+g=0
}
over  global fields, which amounts to the integral representability of $-g$ as
 the binary quadratic form $ax^2+bxy+cy^2$.
Here ``integral" means that we shall restrict the equation over rings of integers
 or more generally $S$-integers.
This problem is an old one which dates back to the theory of Gauss' quadratic forms.
Modern approaches involve the splitting of ideals in quadratic extensions and
 class field theory.
For example, the main theorem of  Cox~\cite{cox} gives a  criterion of the solvability
 of the diophantine equation
\eq{\label{eq_x2ny2}
p=x^2+ny^2
}
 for any positive integer $n$ and prime number $p$.
Using a similar argument, the author and Deng~\cite{rcf} generalized the base field
 of \eqref{eq_x2ny2} to a class of imaginary quadratic fields.
Maciak \cite{maciak2011primes} treated the same problem over rational function fields,
 and gave a similar criterion of the integral  solvability of \eqref{eq_x2ny2}.

Another approach is to use the point of view of arithmetic algebraic geometry.
The integral solvability of an equation  amounts to the existence of integral points
 on the affine scheme defined by it.
Colliot-Th\'el\`ene and Xu \cite{colliot2009brauer} studied
the integral points on homogeneous spaces of semi-simple and simply connected
linear algebraic groups of non-compact type by using the strong approximation theorem
 and the Brauer-Manin obstruction. They also applied the results to the integral
 representation problem of quadratic forms.
Harari~\cite{bmob} showed that the Brauer-Manin obstruction is the only obstruction
 for the existence of integral points of a scheme over the ring of integers of a number
 field, whose generic fiber is a principal homogeneous space (torsor) of a torus.
Although these results are applicable to \eqref{eq_qr}, it can not yield an explicit
 criterion for the integral solvability.

After then  Wei and  Xu  \cite{multi-norm-tori,multip-type}  showed that there exist
 idele groups which are the so-called $\XX$-\emph{admissible subgroups} for determining
 the integral points for multi-norm  tori
 (more generally, groups of multiplicative types),
 and  interpreted the $\XX$-admissible subgroup
 in terms of finite Brauer-Manin obstruction.
In \cite[Section 3]{multi-norm-tori} Wei and Xu also showed how to apply this method
 to binary quadratic diophantine equations over rings of integers of number fields.
As applications,  they gave  some explicit criteria of the solvability of equations of the
 form $x^2\pm dy^2=a$ over $\ZZ$ in \cite[Sections 4 and 5]{multi-norm-tori},
 by constructing  explicit admissible subgroups.
Later  Wei \cite{wei_diophantine} applied the method in \cite{multi-norm-tori}
to  give some additional criteria of the solvability of the diophantine equation
$x^2-dy^2=a$ over $\ZZ$ for some $d$.
He also determined which integers can be written as a sum of two integral squares for some
 of the quadratic fields $\QQ(\sqrt{\pm p})$ (in \cite{wei1}),
 $\QQ(\sqrt{-2p})$ (in \cite{wei2}) and so on.
The author et al. \cite{lv2018intrepqr} also applied the method in \cite{multi-norm-tori}
 to \eqref{eq_qr}  over $\ZZ$ and gave a criterion of the solvability with some additional
 assumptions, by constructing explicit admissible subgroups for \eqref{eq_qr}.

In this text, we treat the equation \eqref{eq_qr} over $k[t]$,
 as an function field  analogue of \cite{lv2018intrepqr}.
We generalize the method in \cite{lv2018intrepqr} to  construct explicit admissible
 subgroups for the equation \eqref{eq_qr}. See  Lemma \ref{lem_lambda_inv_img}.
Specifically, the main result is:
\thmus{
Let $k=\FF_q$ be a finite field of odd characteristic,
 $F=k(t)$ a rational function field , $\fo_F=k[t]$.
Suppose $a,b,c$ and $d$ are elements of $\fo_F$ such that
 $E=F(\sqrt{(b/2)^2-ac})/F$ is an imaginary quadratic extension.
Let $\Kp$ be the  class field   corresponding to $E^\tm\Xp$ and
\eqn{
\XX=\Spec(\fo_F[x,y]/(a(ax^2+bxy+cy^2+g))).
}
Then the equation \eqref{eq_qr}
 is solvable over $\fo_F$  if and only if there exists a local solution
\eqn{
\prod_{\p\in\Omega_F}(x_\p, y_\p)\in\prod_{\p\in\Omega_F}\XX(\fo_{F_\p})
}
such that
\eqn{
\psi_{\Kp/E}(\tilde f_E(\prod_\p(x_\p,y_\p)))=1.
}}
In the above  theorem, ``imaginary" means
 there is a unique place lying over $1/t$,
 $E^\tm\Xp$ (depending on the sign function)
 is an open subgroup  of finite index of the idele group $\II_E$  of $E$,
 $\tilde f_E$ is a map from $\prod_\p\XX(\fo_{F_\p})$ to  $\II_E$
 which is constructed by using the fact that
  the generic fiber of $\XX$ admits the structure of a torsor of a torus, and
$\psi_{\Kp/E}: \II_E\rightarrow \Gal(\Kp/E)$ is the Artin reciprocity map.
The condition $ \psi_{\Kp/E}(\tilde f_E(\prod_\p(x_\p,y_\p)))=1$ is called the Artin condition.
See Sections \ref{sec_nota} and \ref{sec_rff} for details.

In Section \ref{sec_artin_cond}, we introduce  from  \cite{multi-norm-tori}  notations
 and the  general result we mainly use in this text, but in a modified way which focus on our goal.
Then we give our result on the equation  \eqref{eq_qr} over $k[t]$ in Section \ref{sec_rff}.
The results state that
the integral local solvability and the Artin condition
(see Remark \ref{rk_artin_cond})
 completely describe the global integral solvability.
In view of  the result of  Maciak, adding an assumption (see \ref{eq_as_hil_rec}),
 we recover the main theorems
 in  \cite{maciak2011primes} by our result. At last, we ended this text by
 concrete examples  showing  the explicit criteria of the solvability.

\section{Solvability by the Artin Condition}\label{sec_artin_cond}

\subsection{Notations}\label{sec_nota}
Let $k=\FF_q$ and $F/k(t)$ a algebraic  function field with characteristic not $2$, $\fo_F$ integral
 closure of $k[t]$ in $F$, $\Omega_F$ the set of all places in $F$.
Thus $2\in \fo_F^\tm$.
Let $F_\p$ be the completion of $F$ at $\p$ and $\fo_{F_\p}$  the valuation ring
 of $F_\p$ for each $\p\in\Omega_F$.
Denote by $\infty_F\subset\Omega_F$ the set of infinite places of $F$, i.e.,
 places lying over $1/t$.
We also write $\fo_{F_\p}=F_\p$ for $\p\in \infty_F$.
The adele ring (resp.  idele group) of $F$ is denoted by $\bA_F$ (resp. $\II_F$).

Let $a,b,c$ and $g$ be  elements in $\fo_F$ such that $-d=(b/2)^2-ac$ is not a square in $F$.
Let $E=F(\sqrt{-d})$, a quadratic extension of $F$.
Since the characteristic of $F$ is not $2$, the extension $E/F$ is separable.
Let
\eq{\label{eq_XX}
\XX=\Spec(\fo_F[x,y]/(a(ax^2+bxy+cy^2+g)))
}
be the affine scheme defined by the equation
\eq{\label{eq_bqf}
a(ax^2+bxy+cy^2+g)=0
}
 over $\fo_F$.
Since $-d$ is not a square in $F$, we have $a\neq0$.
Then the equation
\eqn{
ax^2+bxy+cy^2+g=0
}
is solvable over $\fo_F$ if and only if $\XX(\fo_F)\neq\emptyset$.

Now we denote
\aln{
\tilde x&= ax+\frac{b}{2}y,	\\
\tilde y&= y,	\\
n&=-ag.
}
Then we can write  \eqref{eq_bqf} as
\eq{\label{eq_bqf_n}
\tilde x^2+d\tilde y^2=n.
}
Denote by $R_{E/F}(\GG_m)$ the Weil restriction of $\GG_{m,E}$ to $F$. Let
\eqn{ \varphi: R_{E/F}(\GG_m)\lra\GG_m }
be the homomorphism of algebraic groups which represents
\eqn{ x\lmpt N_{E/F}(x): (E\otm_FA)^\tm\lra A^\tm}
for any $F$-algebra $A$. Define the torus $T=\ker\varphi$.
Let $X_F$ be the generic fiber of $\XX$.
We may write an  element in $T(A)$ (resp.  $X_F(A)$) as
$u+\sqrt{-d}v\in E\otm_FA$, with $u,v\in A$, $u^2+dv^2 = 1$
 (resp. $\tilde x+\sqrt{-d}\tilde y\in E\otm_FA$ with $x,y\in A$, $\tilde x^2+d\tilde y^2=n$).
Then $X_F$ is naturally a $T$-torsor by the action:
\aln{
T(A)\tm X_F(A) &\lra X_F(A)\\
(u+\sqrt{-d}v, \tilde x+\sqrt{-d}\tilde y) &\lmpt (u+\sqrt{-d}v) (\tilde x+\sqrt{-d}\tilde y).
}
Obviously,  $T$ has an integral model $\TT=\Spec(\fo_F[x,y]/(x^2+dy^2-1))$
 and since $\TT$ is separated over $\fo_F$, we can view $\TT(\fo_{F_\p})$ as a subgroup of $T(F_\p)$.
Note that for $\p\in\infty_F$ we also write $\TT(\fo_{F_\p})$ (resp. $\XX(\fo_{F_\p})$)
 for $T(F_\p)$ (resp. $X_F(F_\p)$).

Denote by $\lambda$ the embedding of $T$ into $R_{E/F}(\GG_m)$.
Clearly $\lambda$ induces a natural injective group homomorphism
\eqn{ \lambda_E: T(\bA_F)\lra\II_E. }

Now we assume that
\eq{\label{eq_nonempty}
X_F(F)\neq\emptyset,
}
i.e. $X_F$ is a trivial $T$-torsor.
Fixing a rational point $P\in X_F(F)$,
for any $F$-algebra $A$, we have an isomorphism
\gan{
\xymatrix{\phi_P: X_F(A)\ar[r]^-{\sim} &T(A)}\\
\qquad\qquad x \lmpt P^{-1}x
}
induced by $P$. That is, $P^{-1}x$ is  the unique element in $T(A)$
 sending $P$ to $x$.
Since  $\XX$ is separated over $\fo_F$,
 we can view $\prod_{\p\in\Omega_F}\XX(\fo_{F_\p})$ as a subset of $X_F(\bA_F)$, and the composition
\eqn{
f_E=\lambda_E\phi_P: \prod_\p\XX(\fo_{F_\p})\lra  \II_E}
 makes sense.
Note that under the previous description for $X_F(A)$ where $A=F$ here,
 we have $P\in E^\tm\subset \II_E$ since it is a rational point.
It follows that we can define the map $\tilde f_E$ to be the composition
\eqn{\xymatrix{
\prod_\p\XX(\fo_{F_\p}) \ar@{->}[r]^-{ f_E}  &\II_E       \ar@{->}[r]^-{\tm P} &\II_E.
}}
It can be seen that the restriction to $\XX(\fo_{F_\p})$ of   $\tilde f_E$ is defined by
\eq{\label{eq_tilde_f_E}
\tilde f_E[(x_\p,y_\p)]=
\cs{
(\tilde x_\p+\sqrt{-d}\tilde y_\p, \tilde x_\p-\sqrt{-d}\tilde y_\p) \in E_\fP\tm E_{\bar\fP}
                                                &\text{if }\p=\fP\bar\fP\text{ splits in }E/F,\\
 \tilde x_\p+\sqrt{-d}\tilde y_\p     \in E_\fP	&\text{otherwise},
}}
where $\fP$ and $\bar\fP$ (resp. $\fP$) are places of $E$ above $\p$
 and $\tilde x_\p=ax_\p+\frac{b}{2}y_\p$, $\tilde y_\p=y_\p$.

Let $\Xi$ be an open subgroup of $\II_E$ such that $E^\tm\Xi$ is of
 finite index.
Let $K_\Xi$ be the class field corresponding to $E^\tm \Xi$ under  class field theory,
such that  the Artin map  gives the isomorphism
\eqn{
\xymatrix{\psi_{K_\Xi/E}: \II_E/E^\tm\Xi\ar[r]^-{\sim} &\Gal(K_\Xi/E).
}}
For any $\prod_{\p\in\Omega_F}(x_\p, y_\p)\in\prod_{\p\in\Omega_F}\XX(\fo_{F_\p})$, noting that $P$ is in $E^\tm$, we have
\eq{\label{eq_tilde_f}
\psi_{K_\Xi/E}(f_E(\prod_\p(x_\p,y_\p)))=1\text{ if and only if }
\psi_{K_\Xi/E}(\tilde f_E(\prod_\p(x_\p,y_\p)))=1.
}

\rk{\label{rk_hasse_min}
The assumption \eqref{eq_nonempty} is easy to check by the Hasse-Minkowski theorem
 on quadratic equations.
In particular, it holds if $\prod_{\p\in\Omega_F}\XX(\fo_{F_\p})\neq\emptyset$,
 in which case, we can pick an $F$-point $P$ of $X_F$ and obtain $\phi_P$.
Note that the map $\tilde f_E$ is independent of $P$.
}

\subsection{A general result}\label{sec_main}
For the integral points of the scheme $\XX$ over $\fo_F$ defined
in \eqref{eq_XX}, we observe that $\TT(\fo_{F_\p})$ acts stably on $\XX(\fo_{F_\p})$ for
 all $\p\in\Omega_F$.
To verify  this, it suffice to show that
 for any $u,v,x,y\in \fo_{F_\p}$,
\eqn{
(u+v\sqrt d)((ax+\frac{b}{2}y)+y\sqrt d)=((ax'+\frac{b}{2}y')+y'\sqrt d)
}
 for some $x',y'\in \fo_{F_\p}$. Indeed, this is the case if we take
\gan{
x'=(u-\frac{bv}{2})x+cvy, \\
y'=uy+v(ax+\frac{b}{2}y).
}
Now we have the following general result, which is
a function field analogue of \cite[Corollary 1.6]{multi-norm-tori}.
\prop{\label{prop_artin_cond}
Let $\Xi$ be an open subgroup of $\II_E$ described as before and suppose that
\eq{\label{eq_lambda_inv_img}
\lambda_E^{-1}(E^\tm\Xi) \subseteq  T(F)\prod_\p \TT(\fo_{F_\p}).
}
Then $\XX(\fo_F)\neq\emptyset$ if and only if there exists a local solution
\eqn{
\prod_{\p\in\Omega_F}(x_\p, y_\p)\in\prod_{\p\in\Omega_F}\XX(\fo_{F_\p})
}
such that
\eq{\label{eq_artin_cond}
\psi_{K_\Xi/E}(\tilde f_E(\prod_\p(x_\p,y_\p)))=1.
}}
\pf{
The proof is the similar to the number field case. See \cite[Corollary 1.6]{multi-norm-tori}.
To imitate the proof, the only nontrivial thing is that
  $\TT(\fo_{F_\p})$ acts stably on $\XX(\fo_{F_\p})$ for  all $\p\in\Omega_F$, as showed above.
}
\rk{\label{rk_artin_cond}
The equation \eqref{eq_artin_cond} is called the \emph{Artin condition} in Wei \cite{wei_diophantine,wei1,wei2}.
If the assumption  in the proposition holds,
 the integral local solvability  and the Artin condition completely describe the global integral solvability.
As a result, in cases where $K_\Xi$ is known it is possible to calculate the Artin condition, and give
  explicit criteria for the solvability.
Actually,  the  idele group $\Xi$ satisfying the
 assumption \eqref{eq_lambda_inv_img} is a variant of the definition of $\XX$-admissible
 subgroup in \cite{multi-norm-tori}.
}

\rk{\label{rk_trivial_artin}
In the case that  $K_\Xi/F$ is abelian, the   Artin condition is trivially true
for any local solution.
Actually, since  $K_\Xi/F$ is abelian, the following diagram commutes:
\eqn{\xymatrix{
\II_E \ar@{->}[rr]^-{\psi_{K_\Xi/E}} \ar@{->}[d]_{N_{E/F}} & &\Gal(K_\Xi/E) \ar@{_(->}[d]\\
\II_F \ar@{->}[rr]^-{\psi_{K_\Xi/F}}                       & &\Gal(K_\Xi/F)
}}
It follows that for any local solution  $\prod_\p(x_\p,y_\p)$,
\eqn{
\psi_{K_\Xi/E}(\tilde f_E(\prod_\p(x_\p,y_\p)))
= \psi_{K_\Xi/F}(N_{E/F}( \tilde f_E(\prod_\p(x_\p,y_\p))))
= \psi_{K_\Xi/F}(n) = 1,
}
where the second equation comes from the definition \eqref{eq_tilde_f_E} of $\tilde f_E$  and
\eqn{
(\tilde x_\p+ \sqrt{-d}\tilde y_\p)(\tilde x_\p- \sqrt{-d}\tilde y_\p)
  = n\text{ in }E_\fP\text{ with }\fP\mid \p,
}
 and the last equality is obtained by the assumption that $n = -ag\in F$.
}

Let $L=\fo_F+\fo_F\sqrt{-d}$ in $E$ and $L_\p=L\otm_{\fo_F}\fo_{F_\p}$
in  $E_\p=E\otm_F F_\p$.
We also write $L_\p=E_\p$ for $\p\in\infty_F$.
Then  $\prod_\p L_\p^\tm$ is an open subgroup of $\II_E$.
Let $S\subseteq \Omega_E$ be a finite set of places of $E$,
$U_\fP\subseteq \fo_{E_\fP}^\tm$ be an open subgroup of $\fo_{E_\fP}^\tm$ for $\fP\in S$  and
\eqn{
W_\fP = \cs{ U_\fP            &\text{ for }\fP\in S,\\
            \fo_{E_\fP}^\tm    &\text{ for }\fP\notin S.
}}
We  define  the open  subgroup of $\II_E$
\eq{\label{eq_Xi_W}
\Xi_W =\left(\prod_{\p\in\Omega_F} L_\p^\tm\right) \bigcap\left( \prod_{\fP\in\Omega_E} W_\fP\right)
      = \prod_\p\left( L_\p^\tm\cap \prod_{\fP\mid\p} W_\fP\right ),
}
and assume  $E^\tm\Xi_W$ is also  of finite index in $\II_E$.

By some additional assumptions, we  prove that
$\Xi=\Xi_W$ satisfies the assumption \eqref{eq_lambda_inv_img} in
Proposition \ref{prop_artin_cond},
that is,
\eqn{
\lambda_E^{-1}(E^\tm\Xi_W) \subseteq  T(F)\prod_\p \TT(\fo_{F_\p}).
}
\lemm{\label{lem_lambda_inv_img}
Let $S$ and $W_\fP$ be as before.
Suppose  for every $u\in\fo_F^\tm$,
the equation
\eqn{
N_{E/F}(\alpha)=u,\quad \alpha\in L^\tm
}
 is solvable or
the equation
\eqn{
N_{E_\p/F_\p}(\alpha)=u,\quad \alpha\in L_\p^\tm \cap \prod_{\fP\mid\p} W_\fP
}
 is not solvable for some place $\p$.
Then  the assumption \eqref{eq_lambda_inv_img} in
{\upshape Proposition \ref{prop_artin_cond}} is true.
}
\pf{
The proof is similar to \cite[Lemma 1]{lv2018intrepqr} but a little different.
Recall that $T=\ker(R_{E/F}(\GG_m)\rightarrow\GG_m)$ and $\TT$ is the group scheme
defined by the equation $x^2+dy^2=1$ over $\fo_F$.
  Therefore we have
\eqn{ T(F)=\{\beta\in E^\tm\mid N_{E/F}(\beta)=1\} }
and
\eqn{ \TT(\fo_{F_\p})=\{\beta\in L_\p^\tm\mid N_{E_\p/F_\p}(\beta)=1\}. }
Suppose   $t\in T(\bA_F)$ such that $\lambda_E(t)\in E^\tm\Xi_W$.
 Write $t=\beta i$ with $\beta\in E^\tm$ and $i\in \Xi_W$.
Since $t\in T(\bA_F)$ we have
\eqn{ N_{E/F}(\beta)N_{E/F}(i)=N_{E/F}(\beta i)=1. }
It follows that
\eqn{ N_{E/F}(i)=N_{E/F}(\beta^{-1})\in F^\tm\cap \prod_\p\fo_{F_\p}^\tm=\fo_F^\tm. }
So we have   $N_{E/F}(i)=u$ for some $u\in\fo_F^\tm$.
Note that $i\in \Xi_W$, and thus at each $\p$ we have
\eqn{
N_{E_\p/F_\p}(i_\p)=u,\quad i_\p=(i_\fP)_{\fP\mid\p}
  \in L_\p^\tm \cap \prod_{\fP\mid\p} W_\fP.
}
Thus the assumption tells us  that
the equation
\eqn{
N_{E/F}(\alpha)=u,\quad \alpha\in L^\tm
}
 is solvable.
Let  $\alpha_0$ be such a solution and let
\aln{
\gamma      &=\beta \alpha_0\\
\text{and }j&=i\alpha_0^{-1}.
}
Then $N_{E/F}(\gamma)=N_{E/F}(j)=1$. Note that $\alpha_0\in L^\tm$,
and we have $\gamma\in T(F)$ and $j\in\prod_\p\TT(\fo_{F_\p})$.
It follows that $t=\beta i=\gamma j\in T(F)\prod_\p\TT(\fo_{F_\p})$.
This finishes the proof.
}
\rk{
In \cite{lv2018intrepqr}, the admissible subgroup $\Xi$ for
 the equation \eqref{eq_qr} is simply chosen to be $\prod_\p L_\p^\tm$,
which is generalized by the above lemma in the function field case, where we intersect
 $\prod_\p L_\p^\tm$ with the  open subgroup $\prod_\fP W_\fP$
 (see \eqref{eq_Xi_W}).
This allows us to deal more difficult base global fields and parameters for \eqref{eq_qr}.
Previous method to do this \cite{multi-norm-tori,wei_diophantine,wei2}
 is to construct an Kummer extension $\Theta/E$ with low degree and
 choose the class group to be
 $E^\tm\prod_\p L_\p^\tm\cap E^\tm N_{\Theta/E}\II_\Theta^\tm$.
}
Using this lemma,  we obtain the following corollary to Proposition
\ref{prop_artin_cond}.
\coru{\label{cor_artin_cond}
Let $\Xi_W$  be defined by \eqref{eq_Xi_W}  and suppose that
$E^\tm\Xi_W$ is of finite index in $\II_E$.
Let $S$ and $W_\fP$ satisfy  the assumption in
Lemma \ref{lem_lambda_inv_img},
that is,  for every $u\in\fo_F^\tm$,
the equation
\eqn{
N_{E/F}(\alpha)=u,\quad \alpha\in L^\tm
}
 is solvable or
the equation
\eqn{
N_{E_\p/F_\p}(\alpha)=u,\quad \alpha\in L_\p^\tm \cap \prod_{\fP\mid\p} W_\fP
}
 is not solvable for some place $\p$.
Then $\XX(\fo_F)\neq\emptyset$ if and only if there exists a local solution
\eqn{
\prod_{\p\in\Omega_F}(x_\p, y_\p)\in\prod_{\p\in\Omega_F}\XX(\fo_{F_\p})
}
such that
\eqn{
\psi_{K_W/E}(\tilde f_E(\prod_\p(x_\p,y_\p)))=1,
}
where $K_W$ is the class field
 corresponding to $E^\tm\Xi_W$.
}

\section{The integral representation of binary quadratic forms over $k[t]$} \label{sec_rff}
Now we consider our focus, the case where $F=k(t)$ and  $k=\FF_q$ is a finite
 field of characteristic $p\neq2$.
Hence we are interested in the diophantine equation
\eqn{
ax^2+bxy+cy^2+g=0
}
 over  $\fo_F=k[t]$.
Suppose that  $-d=(b/2)^2-ac$ is not a square in $F$.
Set $E=F(\sqrt{-d})$ and $L=\fo_F+ \fo_F\sqrt{-d}$ as previous sections.
Let  $\pin=1/t$ be  the place of $k(t)$ at infinity and
 suppose  further that  $E/F$ is ``imaginary", that is,
\eq{\label{eq_as_unique_infty}
\text{there is a unique  place $\Pin$ in $E$ lying over $\pin$.}
}
We briefly introduce sign function here. In a completion $K_\p$ of a global
 function field $K$, A \emph{sign function} with respect to a  uniformizer
 $\pi$ is defined as
\aln{
\sgn: K_\p &\lra \fo_{K_\p}/\p\\
x          &\lmpt c_{r_0},
}
 where $c_{r_0}\neq0$ is the leading coefficient of the Laurent series
 $x=\sum_{r=r_0}^\infty c_r\pi^r$  of $x$ with coefficient in $\fo_{K_\p}/\p$.
An element $x\in K_\p$ is \emph{positive} if $\sgn(x)=1$.
Fix a  sign function of $\Pin$, denoted by $\sgn(\cdot)$ and define the open subgroup
\eq{\label{eq_Ep}
\Ep=\{\alpha\in \EP\mid \sgn(\alpha)=1\}\subseteq \EP
}
 consisting all positive elements, $S=\{\Pin\}$, and $U_\Pin=\Ep$.
Let $\Xp$ be the subgroup $\Xi_W$ defined in \eqref{eq_Xi_W}
 for the chosen $S$ and $U_\Pin$.
\thm{\label{thm_rff}
With the above notations, we have:
\enmt{[\upshape (a)]
\item The open subgroup $E^\tm\Xp$  is of finite index in $\II_E$. \label{it_fin}
\item Let  $\Kp$ be  the  class field   corresponding to $E^\tm\Xp$ and \label{it_sol}
\eqn{
\XX=\Spec(\fo_F[x,y]/(a(ax^2+bxy+cy^2+g))).
}
Then $\XX(\fo_F)\neq\emptyset$ if and only if there exists a local solution
\eqn{
\prod_{\p\in\Omega_F}(x_\p, y_\p)\in\prod_{\p\in\Omega_F}\XX(\fo_{F_\p})
}
such that
\eqn{
\psi_{\Kp/E}(\tilde f_E(\prod_\p(x_\p,y_\p)))=1.
}}}
\pf{
We first show that  $\II_E/E^\tm\Xp$  is finite.
By  the choice of  $S$ and $U_\Pin$ we know that
\eqn{
\Xp= \Ep\tm \prod_{\p\neq\pin} L_\p^\tm
}
Define $(\fo_E)_\p=\fo_E\otm_{\fo_F}\fo_{F_\p}$ for $\p\neq\pin$
in  $E_\p=E\otm_F F_\p$ and
\eqn{
\tXp = \Ep\tm \prod_{\p\neq\pin} (\fo_E)_\p^\tm
     = \Ep\tm \prod_{\fP\neq\Pin} \fo_{E_\fP}^\tm.
}
Since we have a surjection $\tXp/\Xp\lra E^\tm\tXp/E^\tm\Xp$
and
\eqn{
\tXp/\Xp=\prod_{\p\mid [\fo_E:L]}(\fo_E)_\p^\tm/L_\p^\tm
}
is finite, we know that $E^\tm\tXp/E^\tm\Xp$ is finite.
Therefore we only need to show that $\II_E/E^\tm\tXp$  is  finite.

Define $\II_E^+=\II_E\cap \Ep$, the subgroup of $\II_E$ consisting
 elements  whose   projection to $\EP$ is in $\Ep$.
Let $E^+=\II_E^+\cap E^\tm$.
Then naturally we have an isomorphism $\II_E^+/E^+\tXp\cong Cl^+(\fo_E)$,
 where $Cl^+(\fo_E)$ is the
 \emph{narrow class group} with respect to  $(\Pin,\sgn)$, which is finite
 (c.f.~\cite[p.~200]{goss1996basicff}).
Since $\II_E=E^\tm\II_E^+$ by
 weak approximation theorem (c.f. \cite[Chapter II.6]{cassels-nt}),
 we have $\II_E^+/E^+\tXp\cong \II_E/E^\tm\tXp$
and thus $\II_E/E^\tm\tXp\cong Cl^+(\fo_E)$ is finite.
This completes the proof for \eqref{it_fin}.

For the assertion \eqref{it_sol} we  apply Corollary \ref{cor_artin_cond}.
For any $u\neq1$ in $\fo_F^\tm=k^\tm$.
Let $\p=\pin$ so $L_\p^\tm \cap \prod_{\fP\mid\p} W_\fP=\Ep$,
since $\Pin$ is the only place above $\pin$.
Assume  that there is $\alpha \in \Ep$  such that
 $u=N_{E_\p/F_\p}(\alpha)=\alpha\bar\alpha$.
By the definition of $\Ep$ \eqref{eq_Ep} we know that $\sgn(\alpha)=\sgn(\bar\alpha)=1$.
It follows that  $u=\sgn(u)=\sgn(\alpha)\sgn(\bar\alpha)=1$,
 which is a contradiction and  shows that the equation
\eqn{
N_{E_\p/F_\p}(\alpha)=u,\quad \alpha\in L_\p^\tm \cap \prod_{\fP\mid\p} W_\fP
}
 is not solvable for  $\p=\pin$.
Thus \eqref{it_sol} follows form Corollary \eqref{cor_artin_cond}.
}

In contrast to the class field $\Kp$ corresponding to $\Xp$, we denote $\KP$
 the Hilbert class field of $E$ with respect to  $\Pin$, i.e.
 $\KP$ is the class field corresponding to the open subgroup
 $E^\tm\XP$ of finite index in $\II_E$,  where
\eqn{
\XP = \EP\tm \prod_{\fP\neq\Pin} \fo_{E_\fP}^\tm,
}
which we will use later and
 basically we have $\Kp\supseteq\KP$ since $E^\tm\Xp\subseteq E^\tm\XP$.

We use  the above theorem to  derive a  result  similar to the
 main theorems  of Maciak
 \cite{maciak2011primes} considering the equation  $l=x^2+Dy^2$.
First we need some notations and facts in \cite{maciak2011primes}.
Recall that we set $F=k(t)$, $k=\FF_q$ and hence  $\fo_F=k[t]$.
Let $D\in k[t]$  be square free with positive degree and
 $l\nmid D$ an irreducible  element of $k[t]$ and
 we consider the equation $l=x^2+Dy^2$ over $k[t]$.
For this equation, we have $a=1$, $b=0$, $c=D$, $g=-l$, $-d=(b/2)^2-ac=-D$
 and  $E = \QQ(\sqrt{-d})$. Thus $\tilde x=x$, $\tilde y=y$ and $n=l$.
Remember that $D$ is square free with positive degree,
 we know that $\fo_F+ \fo_F\sqrt{-D}=\fo_E$.
Recall that \cite[Proposition 14.6]{rosen2013number}
 $\pin$ ramifies,  splits, or is innert in $E/F$ if $\deg D$ is odd, $\deg D$
 is even and $\lc(-D)\in k^{\tm2}$, or $\deg D$ is even and $\lc(-D)\not\in k^{\tm2}$,
 respectively, where $\lc$ means the leading coefficient.
Suppose that $\deg D$ is odd or $\lc(-D)\not\in k^{\tm2}$,
 which is to say  the assumption
 \eqref{eq_as_unique_infty} holds.
Then a necessary condition for $l=x^2+Dy^2$ being solvable over $k[t]$ is
\eq{\label{eq_as_deg_d_even}
\text{$\deg l$ is even if $\deg D$ is.}
}
We will always assume this.
Let  $d_\infty$ be the relative degree of $\Pin\mid\pin$ and
 define $\deg^*l=\frac{\deg l}{d_\infty}$ as in \cite[Section 4]{maciak2011primes},
 which is an positive integer by the above assumption.
Let $g$ be the genus of $E$. Following \cite{maciak2011primes}
 we fix $\sgn$  with respect to the uniformizer $t^g/\sqrt{-D}$
 (i.e., $\sgn(t^g/\sqrt{-D})=1$)
 in the following
\thm{\label{thm_h+}
Let $k$,  $D$ and $l$ in $k[t]$ be as before such that \eqref{eq_as_deg_d_even}
 holds. Then  we have
\enmt{[\upshape (a)]
\item if  $\sgn(l)(-1)^{\deg^*l}\in k^{\tm2}$, \label{it_sq}
 then $l=x^2+Dy^2$ is solvable over $k[t]$ if and only if
 $\fracn{l}{r}=1$ for each monic irreducible factor $r\mid D$ and
 $l$ splits completely in $\Kp$;
\item if  $\sgn(l)(-1)^{\deg^*l}\not\in k^{\tm2}$, \label{it_nsq}
 then $l=x^2+Dy^2$ is solvable over $k[t]$ if and only if
 $\fracn{l}{r}=1$ for each monic irreducible factor $r\mid D$,
 $l$ splits completely in $\KP$ and
 the relative degree of $l$ in $\Kp$ is $2$.
}}
\pf{
In line with  Theorem \ref{thm_rff},
 recall that $F=k(t)$, $E=F(\sqrt{-D})$,
 $L=\fo_F+ \fo_F\sqrt{-D}=\fo_E$ and  $\Kp$ (resp. $\KP$) is the class field
 corresponding to $\Xp$ (resp. $\XP$).
We know by \eqref{eq_tilde_f_E} that
\eq{\label{thm_recover.eq_tilde_f_E}
\tilde f_E[(x_\p,y_\p)]=
\cs{
(x_\p+\sqrt{-D}y_\p, x_\p-\sqrt{-D}y_\p) &\text{if }\p\text{ splits in }E/F,\\
 x_\p+\sqrt{-D}y_\p                      &\text{otherwise}.
}}
Then by Theorem \ref{thm_rff}, the equation
$l=x^2+Dy^2$ is solvable over $k[t]$  if and only if there exists a
local solution
\eqn{
\prod_{\p\in\Omega_F}(x_\p, y_\p)\in\prod_{\p\in\Omega_F}\XX(\fo_{F_\p})
}
such that
\eqn{
\psi_{\Kp/E}(\tilde f_E(\prod_\p(x_\p,y_\p)))=1.
}

Next we verify these conditions in details.
By  a simple  calculation we know the local condition  \eqn{
\prod_{\p}\XX(\fo_{F_\p})\neq\emptyset
} is equivalent to
\enmt{[(I)]
\item \label{thm_h+.eq_local_infty}
\hspace{2.5cm} \text{$\lc(l)\lc(D)^{\deg l}\in k^{\tm2}$ if $\deg D$ is odd},
}
\al{
&\text{$l$ splits completely in $E$},\label{thm_h+.eq_local_l}\\
&\fracn{l}{r}=1,\text{ for each monic irreducible factor $r\mid D$.}\label{thm_h+.eq_local_r}
}

For the Artin condition, let $\prod_\p(x_\p,y_\p)\in\prod_\p\XX(\fo_{F_\p})$
be a local solution. Then
\eq{\label{thm_recover.eq_local_decomp}
(x_\p+\sqrt{-D}y_\p)(x_\p-\sqrt{-D}y_\p) = l
  \text{ in }E_\fP\text{ with }\fP\mid\p.
}
Let $\fl=l\fo_F$.
Thus for all $\p\nmid \fl\pin$, $\tilde f_E[(x_\p,y_\p)]\in L_\p^\tm$ by
\eqref{thm_recover.eq_tilde_f_E} and  \eqref{thm_recover.eq_local_decomp}.
It follows that
\eqn{
\psi_{\Kp/E}(\tilde f_E[(x_\p,y_\p)])=1\quad\text{for all $\p\nmid \fl\pin$},
}
where $\tilde f_E[(x_\p,y_\p)]$ is  regarded as an element in $\II_E$ such that
 the component above $\p$ is given by the value of $\tilde f_E[(x_\p,y_\p)]$ and $1$ otherwise.
For $\p=\fl$,   by the local condition we already know that
$\fl$ splits completely in $E/F$.
Hence \eqref{thm_recover.eq_local_decomp} tells us that one of
$v_\fl(x_\fl\pm\sqrt{-D}y_\fl)$  is $1$ and the other $0$.
Suppose  $v_\fl(x_\fl+\sqrt{-D}y_\fl)=1$ and
let $\fl=\fL\bar\fL$ in $E$.
Note that $L=\fo_E$ and
 $L_\fl^\tm = (\fo_E)_\fl^\tm = \fo_{E_\fL}^\tm\tm \fo_{E_{\bar\fL}}^\tm$,
 so both $\fL$ and $\bar\fL$ are unramified in $\Kp/E$.
It follows that
\eqn{
\s_\fL=\psi_{\Kp/E}(\tilde f_E[(x_\fl,y_\fl)])=\psi_{\Kp/E}(l_\fL)\in\Gal(\Kp/E)
}
where $l_\fL$ is in $\II_E$ such that its $\fL$ component is $l$ and the others
 $1$, and $\s_\fL$ denotes
 the Frobenius automorphism  of $\fL$ in $\Kp/E$.

For $\p=\pin$,  we have   $l=\alpha\bar\alpha$ where
\eqn{
\alpha=\tilde f_E[(x_\pin,y_\pin)]= x_\pin+\sqrt{-D}y_\pin.
}
Then   $\sgn(\bar\alpha)=(-1)^{\deg^*l}\sgn(\alpha)$
 (\cite[Proposition 4.3]{maciak2011primes}).
Thus
\eqn{
\sgn(\alpha)=\pm\sqrt{\sgn(l)(-1)^{\deg^*l}}.
}
Let
\eqn{
 \s_\Pin =\psi_{\Kp/E}(\tilde f_E[(x_\pin,y_\pin)])\in\Gal(\Kp/E).
}
So we obtain that  the Artin condition
$\psi_{\Kp/E}(\tilde f_E(\prod_\p(x_\p,y_\p)))=1$ is equivalent to
\eq{\label{thm_h+.eq_artin}
\s_\fL\s_\Pin=1.
}

Note that since $L=\fo_E$, we have
\eqn{
\Xp= \Ep\tm \prod_{\p\neq\pin} L_\p^\tm=\Ep\tm \prod_{\fP\neq\Pin} \fo_{E_\fP}^\tm.
}
It follows that for any $\beta\in \EP$, $\beta\in  E^\tm\Xp$ if and only if
 $\sgn(\beta)\in \fo_E^\tm=k^\tm$.
We also note that  $\sgn(l)(-1)^{\deg^*l}\in k^\tm$ since
 $\sgn(l)=\lc(l)\lc(-D)^{-\deg^*l}$ by (1) of \cite[p. 230]{maciak2011primes}.
Using these facts, we distinguish two cases:
\enmt{[(i)]
\item $\sgn(l)(-1)^{\deg^*l}\in k^{\tm2}$. Then
 $\sgn(\alpha)=\pm\sqrt{\sgn(l)(-1)^{\deg^*l}}\in k^\tm$
 and thus $\alpha_\Pin\in E^\tm\Xp$.
 It follows that  $\s_\Pin=1$.
\item $\sgn(l)(-1)^{\deg^*l}\not\in k^{\tm2}$.
 But we already know that $\sgn(l)(-1)^{\deg^*l}\in k^\tm$.
 It follow that $\sgn(\alpha)\not\in k^\tm$ and $\sgn(\alpha^2)\in k^\tm$,
 which is to say that
 $\alpha_\Pin$ is of order $2$ in $\II_E/E^\tm\Xp\cong \Gal(\Kp/E)$  and
 hence $\s_\Pin$  is an element of  order $2$.
}

At this time we know that $l=x^2+Dy^2$ is solvable over $k[t]$  if and only if
 \eqref{thm_h+.eq_local_infty},
 \eqref{thm_h+.eq_local_l},
 \eqref{thm_h+.eq_local_r} hold and there is a local solution
\eqn{
\prod_{\p\in\Omega_F}(x_\p, y_\p)\in\prod_{\p\in\Omega_F}\XX(\fo_{F_\p})
}
such that the corresponding
 \eqref{thm_h+.eq_artin} holds. We will use this equivalence in the sequel.

We first consider \eqref{it_sq}. Suppose that $\sgn(l)(-1)^{\deg^*l}\in k^{\tm2}$.
Then $\s_\Pin=1$.
If $l=x^2+Dy^2$ is solvable over $k[t]$,  then by  \eqref{thm_h+.eq_local_l}
 $l$ splits in $E/F$.
Moreover,  \eqref{thm_h+.eq_artin} implies $\s_\fL=\s_\fL\s_\Pin=1$,
 i.e. $\fL$ splits completely in $\Kp/E$.
Thus  $l$ splits completely in $\Kp$ and
 $\fracn{l}{r}=1$ for each monic irreducible factor $r\mid D$
 (which is \eqref{thm_h+.eq_local_r}, and the same for the sequel).
Conversely, if  $l$ splits completely in $\Kp$ and
 $\fracn{l}{r}=1$ for each monic irreducible factor $r\mid D$,
 then $l$ also splits in $E$, i.e. \eqref{thm_h+.eq_local_l} holds.
Also we have $\s_\fL\s_\Pin=\s_\fL=1$ so $\eqref{thm_h+.eq_artin}$ holds.
At last note that $\sgn(l)(-1)^{\deg^*l}\in k^{\tm2}$. Hence if $\deg D$ is odd,
 $\deg^*l=\deg l$ and
\eqn{
\sgn(l)(-1)^{\deg^*l}=\lc(l)\lc(D)^{-\deg^*l}=\lc(l)\lc(D)^{-\deg l}\in k^{\tm2}.
}
Thus \eqref{thm_h+.eq_local_infty} holds.
This completes the proof for \eqref{it_sq}.

To show \eqref{it_nsq}, suppose that $\sgn(l)(-1)^{\deg^*l}\not\in k^{\tm2}$.
Then $\s_\Pin$ is of order $2$.
However, since $\tilde f_E[(x_\pin,y_\pin)]\in\XP$,
\eq{\label{thm_h+.eq_spin_res}
\s_\Pin|_\KP=\psi_{\KP/E}(\tilde f_E[(x_\pin,y_\pin)])=1.
}
If $l=x^2+Dy^2$ is solvable over $k[t]$,  then by  \eqref{thm_h+.eq_local_l}
 $l$ splits in $E/F$.
By \eqref{thm_h+.eq_artin} we know that  $\s_\fL=\s_\Pin^{-1}$ is of order $2$.
Also, by \eqref{thm_h+.eq_spin_res},
$\s_\fL|_\KP=(\s_\Pin|_\KP)^{-1}=1$.
So we have  $l$ splits completely in $\KP$ and
 the relative degree of $l$ in $\Kp$ is $2$.
Conversely, if  $l$ splits completely in $\KP$ and
 the relative degree of $l$ in $\Kp$ is $2$, then
 \eqref{thm_h+.eq_local_l} holds for the same reason as in the prove for \eqref{it_sq}.
In addition,  if $\deg D$ is odd, \eqn{
\sgn(l)(-1)^{\deg^*l}=\lc(l)\lc(D)^{-\deg^*l}=\lc(l)\lc(D)^{-\deg l}\not\in k^{\tm2},
}
 which is impossible. Thus   \eqref{thm_h+.eq_local_infty} trivially holds.
At last it suffices to show $\eqref{thm_h+.eq_artin}$ for a local solution.
Actually, since $l$ splits completely in $\KP$, we have $\s_\fL|_\KP=1$.
It follows that $\s_\fL\in\Gal(\Kp/\KP)$.
Also  $\s_\fL$ is of order $2$ since
 the relative degree of $l$ in $\Kp$ is $2$.
On the other hand,
recall that $\s_\Pin$ is of order $2$ in this case and
we have  $\s_\Pin\in\Gal(\Kp/\KP)$
by \eqref{thm_h+.eq_spin_res}.
Note that  $\Gal(\Kp/\KP)$ is cyclic
 (c.f. \cite[Proposition 7.4.10]{goss1996basicff}).
Therefore $\Gal(\Kp/\KP)$ has a unique subgroup $\{\pm1\}$ of  order $2$
 and $\s_\fL=\s_\Pin=-1\in \{\pm1\}$.
It follow that  $\s_\fL\s_\Pin=1$,
 i.e. that $\eqref{thm_h+.eq_artin}$ holds.
The proof for \eqref{it_nsq} is finished.
}
\rk{
If we assume that
\eq{\label{eq_as_hil_rec}
\text{$\deg D$ is odd or $D$ contains no odd degree irreducible factor},
}
then  \cite[Theorems 4.4 and 4.6]{maciak2011primes}
 is a special form of  the above theorem in the case  $\sgn(l)\in k^{\tm2}$.
To see this it suffice to show that under the assumptions \eqref{eq_as_hil_rec} and
 $\sgn(l)\in k^{\tm2}$,
 the local condition \eqref{thm_h+.eq_local_r} is redundant.
We will use a similar argument as in $F=\QQ$ case \cite[Corollary 4.2]{multi-norm-tori}.
Suppose that the solvable conditions in Theorem \ref{thm_h+} hold.
For monic irreducible $r\mid D$, assuming  \eqref{eq_as_hil_rec} and
 $\sgn(l)\in k^{\tm2}$,
 quadratic reciprocity law implies that
 there exists $u\in \fo_F^\tm$ such that $\fracn{l}{r}=\fracn{ur}{l}$,
 and that one of $ur$ and $-D/(ur)$ have even degree and leading coefficient in $k^{\tm2}$.
 Let $r^*=ur$. We see   that $\sqrt{r^*}\in\KP$ if and only if
\eqn{
\psi_E(i)(\sqrt{r^*})=\sqrt{r^*}\quad\text{for all } i\in \XP
}
under the Artin map $\psi_E$ of $E$,
which is equivalent to the product of  quadratic Hilbert symbols
\eqn{
\prod_\fP \fracn{r^*, i_\fP}{\fP}=1\quad\text{for all $i=(i_\fP)_\fP\in\XP$}.
}
Clearly $E(\sqrt{r^*})/E$ is  unramified  at $\fP\nmid \fr\pin$ where $\fr = r^*\fo_F$.
Since one of $r^*$ and $-D/r^*$ has even degree and leading coefficient in $k^{\tm2}$,
 $\pin$ splits in one of  $F(\sqrt{r^*})$  and $F(\sqrt{-D/r^*})$.
It follows that $\Pin$  splits  in $E(\sqrt{r^*})/E$ and then
 $\fracn{r^*, i_\fP}{\fP}=1$ for all $\fP\nmid \fr$.
Since $i\in \XP = \EP\tm \prod_{\p\neq\pin} (\fo_E)_\p^\tm$ and
 $(\fo_E)_\p=\fo_E\otm_{\fo_F}\fo_{F_\p}
  =L\otm_{\fo_F}\fo_{F_\p}$,
 there exist $a_\p,b_\p\in \fo_{F_\p}$ for each $\p\neq\pin$ such that
\eqn{
\cs{
(i_\fP,i_{\bar\fP}) = (a_\p+\sqrt{-D}b_\p, a_\p-\sqrt{-D}b_\p)
                           &\text{if }\p=\fP\bar\fP\text{ splits in }E/F,\\
i_\fP               =  a_\p+\sqrt{-D}b_\p
                           &\text{otherwise}.
}}
It follows that
\eqn{
\prod_{\fP\mid \fr} \fracn{r^*, i_\fP}{\fP}
  =\prod_{\p\mid\fr}\fracn{r^*, a_\p^2+Db_\p^2}{\p}=1,
}
where the last equality comes from \cite[Ch.V~(3.4)~Proposition]{neukirch_alnt}.
Thus we have  $\sqrt{r^*}\in\KP$ and then $l$ splits in $F(\sqrt{r^*})$
 since it does in $\KP$,
 which is to say $\fracn{l}{r}=\fracn{r^*}{l}=1$
 for each monic irreducible factor $r\mid D$.
This ensures the local condition \eqref{thm_h+.eq_local_r}.
}

We now give two examples where the explicit  criteria  are
 obtained using Theorem \ref{thm_rff}.
\eg{
Let $k=\FF_3$ and $g\in k[t]$,  write
\eqn{
g=u\tm (t-1)^{s_1}\tm (t^2-t-1)^{s_2}\tm \prod_{j=1}^r p_j^{m_j},
}
 where $u\in k^\tm,s_1,s_2,r\ge0,m_j\ge1,p_1,p_2,\dots,p_r \neq t-1, t^2-t-1$ are distinct
 monic irreducible polynomial in $k[t]$.
Then the diophantine equation
\eq{\label{eq_eg}
-x^2+txy-(t^3-t^2+1)y^2+g=0
}
 is solvable over  $k[t]$ if and only if
\enmt{[\upshape (1)]
\item $\fracn{g\tm p^{-v_p(g)}}{p} = (-1)^{v_p(g)}$, for $p=t-1$ or $t^2-t-1$, \label{it_l1}
\item $\fracn{-(t-1)(t^2-t-1)}{p}=1$, for $p\nmid (t-1)(t^2-t-1)$ with \label{it_l2}
 odd $v_p(g)$.
}}
\pf{
In this example, we have $a=-1$, $b=t$, $c=-(t^3-t+1)$ and $-d=(b/2)^2-ac=-(t-1)(t^2-t-1)$.
Since  $\deg d=3$  odd, the assumption \eqref{eq_as_unique_infty} holds.
In order to apply   Theorem \ref{thm_rff},
let  $F=k(t)$, $E=F(\sqrt{-d})$, $L=\fo_F+ \fo_F\sqrt{-d}$ and $\Kp$ is the class field
 corresponding to $\Xp$, with respect to a fixed sign function.
Since $d$ is square free and $\deg d$ is odd,
 we know that  $L=\fo_E$ and $\Kp=\KP$ is the Hilbert class field
 corresponding to
\eqn{
\XP = \EP\tm \prod_{\fP\neq\Pin} \fo_{E_\fP}^\tm
}
(c.f. \cite[Proposition 7.4.10]{goss1996basicff}).
Moreover, if the $j$-invariant of the Drinfeld
 $\fo_E$-module corresponding to the lattice $\fo_E$, the Hilbert class
 field $\KP$ is generated by $j$ over $E$. For proofs,
 see \cite{gekeler1983arithmetik}.
A concise introduction to the theory of Drifeld modules can be found in
 Hayes \cite{hayes1992brief}.
Using Magma calculator \cite{magma} we have $\Kp=\KP=E(\sqrt{t^2-t-1})$.
Thus $\Kp/F$ is abelian, and the Artin condition is trivially true by
Remark \ref{rk_trivial_artin}.
Then by Theorem \ref{thm_rff}, the equation \eqref{eq_eg}
 is solvable over $k[t]$  if and only if it is locally solvable.
By  a simple  calculation we know the local condition
 is equivalent to \eqref{it_l1} and \eqref{it_l2}.
The proof is complete.
}

\eg{
Let $k=\FF_3$ and $g\in k[t]$,  write
\eqn{
g=u\tm q_1^{s_1}\tm q_2^{s_2}\tm \prod_{j=1}^r p_j^{m_j},
}
 where $u\in k^\tm,s_1,s_2,r\ge0,m_j\ge1$, $q_1=t-1,q_2=t^2+t-1$,
 and $p_1,p_2,\dots,p_r \neq q_1,q_2$ are distinct
 monic irreducible polynomial in $k[t]$.
Let $-d=-q_1q_2$, $\theta(X)=X^4-(t^2-t)X^2-t^3+1\in k[t][X]$ and
\aln{
D_1&=\{p=p_1,\dots,p_r\mid \fracn{-d}{p}=1\text{ and }\theta(X)\mod p
 \text{ factors into } \\
 &\qquad\text{ two irreducible polynomials}\},\\
D_2&=\{p=p_1,\dots,p_r\mid \fracn{-d}{p}=1\text{ and }\theta(X)\mod p
 \text{ is irreducible}\}.
}
Then the diophantine equation
\eq{\label{eq_eg2}
(t-1)x^2+(t^2+t-1)y^2+g=0
}
 is solvable over  $k[t]$ if and only if
\enmt{[\upshape (1)]
\item $\fracn{g\tm p^{-v_p(g)}}{p} = (-1)^{\deg(p)}$, for $q_1$ \label{eg2.it_l1}
 or $q_2$,
\item $\fracn{-d}{p}=1$, for $p\nmid d$ \label{eg2.it_l2}
 with odd $v_p(g)$,
\item  and \label{eg2.it_artin}
\aln{
\text{$D_2=\emptyset$ and }&
\sum_{p\in\{q_1,q_2\}\cup D_1} v_p(g)\equiv1\pmod2,\\
\text{or $D_2\neq\emptyset$ and }&
\sum_{p\in D_2} v_p(g)\equiv0\pmod2.
}
}}
\pf{
We have $a=q_1$, $b=0$, $c=q_2$  and $-d=-q_1q_2$.
Since $-d$ is square free and $\deg d$ is odd, we have
 $F=k(t)$, $E=F(\sqrt{-d})$, $L=\fo_F+ \fo_F\sqrt{-d}$ and $\Kp=\KP$ is the
 Hilbert class field as in the previous example.
Using Magma calculator \cite{magma} we have $\Kp=\KP=E[X]/\theta(X)$.
Then by Theorem \ref{thm_rff}, the equation \eqref{eq_eg2}
 is solvable over $k[t]$  if and only if it is locally solvable and the
 Artin condition holds.
It is easy to see that the local condition  is equivalent to
\eqref{eg2.it_l1} and \eqref{eg2.it_l2}.

For the Artin condition, first we know that the discriminant of $\theta$ is $-q_1^4q_2$.
Note that \eqn{
\Gal(\Kp/E)\cong\langle\sqrt{-1}\rangle
}
 is cyclic of order $4$.
This can be shown as follows. Let $p_0=t^2+1$. It is irreducible over $k[x]$
 and we have $\fracn{-d}{p_0}=1$.  Thus $\p_0\fo_E=\fP_0\bar\fP_0$  splits.
One can check that $\theta(X)$ is irreducible over
 $\fo_E/\fP_0\cong\fo_F/p_0\fo_F$, it follows that the Frobenius automorphism
 $\fracn{\Kp/E}{\fP_0}$ has order $4$. Hence $\Gal(\Kp/E)$ is cyclic of order $4$.
Next we only give a sketch of the calculation of Artin condition since it is
 very similar to \cite[Example 1]{lv2018intrepqr} over $F=\QQ$.
In the sequel we identify each finite place $\p$ of $F$ as the unique monic
 irreducible polynomial in $k[x]$ that generates it. We also write $\infty$ for $\pin$.

Let $\XX=\Spec(\fo_F[x,y]/((t-1)((t-1)x^2+(t^2+t-1)y^2+g)))$,
 and $(x_p,y_p)_p\in\prod_p\XX(\fo_{F_p})$.
Note that $\tilde x=q_1x$, $\tilde y=y$ and $n=-q_1g$.
\enmt{[(i)]
\item If $p=q_1\fo_F$, then
 $\psi_{\Kp/E}(\tilde f_E[(x_p,y_p)])=(-1)^{v_p(g)+1}$.
\item If $p=q_2\fo_F$, then
 $\psi_{\Kp/E}(\tilde f_E[(x_p,y_p)])=(-1)^{v_p(g)}$.
\item If $\fracn{-d}{p}=1$ and $\theta(X)\mod p$ splits into linear factors, then
 $\psi_{\Kp/E}(\tilde f_E[(x_p,y_p)])=1$.
\item If $\fracn{-d}{p}=1$ and $\theta(X)\mod p$ splits into two irreducible factors,
 then \eqn{
\psi_{\Kp/E}(\tilde f_E[(x_p,y_p)])=(-1)^{v_p(g)}.
}
\item If $\fracn{-d}{p}=1$ and $\theta(X)\mod p$ is irreducible,
 then \eqn{
\psi_{\Kp/E}(\tilde f_E[(x_p,y_p)])=\pm(\sqrt{-1})^{v_p(g)}
}
 with the sign chosen freely.
\item If $\fracn{-d}{p}=-1$, then $\psi_{\Kp/E}(\tilde f_E[(x_p,y_p)])=1$.
\item If $p=\infty$, then $\psi_{\Kp/E}(\tilde f_E[(x_p,y_p)])=1$ since $\infty$ splits
 completely in $\Kp/E$.
}
Thus the Artin condition $\psi_{\Kp/E}(\tilde f_E(\prod_\p(x_\p,y_\p)))=1$
 is exactly \eqref{eg2.it_artin}.
The proof is complete.
}

\section*{Acknowledgment}
The author would like to thank Jianing Li and Yupeng Jiang for  helpful discussions
 and the referees for valuable suggestions.

This work was supported by
 National Natural Science Foundation of China (Grant No. 11701552).

\bibliography{\bibfilename}

\providecommand{\bysame}{\leavevmode\hbox to3em{\hrulefill}\thinspace}
\providecommand{\MR}{\relax\ifhmode\unskip\space\fi MR }
\providecommand{\MRhref}[2]{%
  \href{http://www.ams.org/mathscinet-getitem?mr=#1}{#2}
}
\providecommand{\href}[2]{#2}
\begin{thebibliography}{10}

\bibitem{magma}
Computational Algebra~Group at~the University~of Sydney,
  \emph{\upshape{Magma}}, V2.23-9, A computational algebra system, available at
  \url{http://magma.maths.usyd.edu.au/magma/}, 2018.

\bibitem{cassels-nt}
J.~W.~S Cassels and A~Fr{\"o}hlich (eds.), \emph{Algebraic number theory},
  Academic Press, 1967.

\bibitem{colliot2009brauer}
Jean-Louis Colliot-Th{\'e}l{\`e}ne and Fei Xu, \emph{Brauer--manin obstruction
  for integral points of homogeneous spaces and representation by integral
  quadratic forms}, Compositio Mathematica \textbf{145} (2009), no.~2,
  309--363.

\bibitem{cox}
David~A. Cox, \emph{Primes of the form {$x^2+ ny^2$}: Fermat, class field
  theory, and complex multiplication}, John Wiley \& Sons, 1989.

\bibitem{gekeler1983arithmetik}
Ernst-Ulrich Gekeler, \emph{Zur arithmetik von {D}rinfeld-moduln},
  Mathematische Annalen \textbf{262} (1983), no.~2, 167--182.

\bibitem{goss1996basicff}
David Goss, \emph{Basic structures of function field arithmetic},
  Springer-Verlag, 1996.

\bibitem{bmob}
David Harari, \emph{Le d{\'e}faut d'approximation forte pour les groupes
  alg{\'e}briques commutatifs}, Algebra \& Number Theory \textbf{2} (2008),
  no.~5, 595--611.

\bibitem{hayes1992brief}
David~R Hayes, \emph{A brief introduction to {D}rinfeld modules}, The
  arithmetic of function fields, vol.~2, de Gruyter Berlin, 1992, pp.~1--32.

\bibitem{rcf}
Chang Lv and Yingpu Deng, \emph{On orders in number fields: {P}icard groups,
  ring class fields and applications}, Science China Mathematics \textbf{58}
  (2015), no.~8, 1627--1638.

\bibitem{lv2018intrepqr}
Chang Lv, Junchao Shentu, and Yingpu Deng, \emph{On the integral representation
  of binary quadratic forms and the {A}rtin condition}, Tokyo J. Math.
  \textbf{41} (2018), no.~2, 371--384. \MR{3908800}

\bibitem{maciak2011primes}
Piotr Maciak, \emph{Primes of the form {$X^2+ nY^2$} in function fields and
  {D}rinfeld modules}, J. Ramanujan Math. Soc. \textbf{26} (2011), no.~2,
  219--235.

\bibitem{neukirch_alnt}
J\"urgen Neukirch, \emph{Algebraic number theory}, Springer, 1999.

\bibitem{rosen2013number}
Michael Rosen, \emph{Number theory in function fields}, Graduate Studies in
  Mathematics, vol. 210, Springer Science \& Business Media, 2013.

\bibitem{wei1}
Dasheng Wei, \emph{On the sum of two integral squares in quadratic fields
  {$\mathbb Q(\sqrt{\pm p})$}}, Acta Arith. \textbf{147} (2011), no.~3,
  253--260.

\bibitem{wei_diophantine}
\bysame, \emph{On the {Diophantine} equation {$x^2-Dy^2=n$}}, Science China
  Mathematics \textbf{56} (2013), no.~2, 227--238.

\bibitem{wei2}
\bysame, \emph{On the sum of two integral squares in the imaginary quadratic
  field {$\mathbb Q(\sqrt{-2p})$}}, Science China Mathematics \textbf{57}
  (2014), no.~1, 49--60.

\bibitem{multi-norm-tori}
Dasheng Wei and Fei Xu, \emph{Integral points for multi-norm tori}, Proceedings
  of the London Mathematical Society \textbf{104} (2012), no.~5, 1019--1044.

\bibitem{multip-type}
\bysame, \emph{Integral points for groups of multiplicative type}, Advances in
  Mathematics \textbf{232} (2013), no.~1, 36--56.

\end{thebibliography}
\bibliographystyle{amsplain}
\end{document}